\documentclass[12pt,reqno]{amsart}
\usepackage{amsmath,amssymb,geometry,color,tikz-cd}
\usepackage[backref,pagebackref,pdftex,hyperindex]{hyperref}
\usepackage[alphabetic]{amsrefs}
\usepackage{amssymb}
\usepackage{bbm}
\usepackage{enumitem}
\usepackage{upgreek}

\newtheorem{proposition}{Proposition}[section]
\newtheorem{theorem}[proposition]{Theorem}
\newtheorem{lemma}[proposition]{Lemma}

\theoremstyle{definition}
\newtheorem{remark}[proposition]{Remark}
\newtheorem{definition}[proposition]{Definition}
\newtheorem{example}[proposition]{Example}

\title{GIT-thresholds for K-moduli spaces of Hypersurfaces}

\author{Chuhui Gao}
\address{School of Mathematical Sciences, Xiamen University, Siming South Road 422, Xiamen, Fujian 361005, China.}
\email{19020241153796@stu.xmu.edu.cn}

\author{Chuyu Zhou}
\address{School of Mathematical Sciences, Xiamen University, Siming South Road 422, Xiamen, Fujian 361005, China.}
\email{chuyuzhou@xmu.edu.cn, chuyuzhou1@gmail.com}

\date{} 

\newcommand{\ord}{{\rm {ord}}}

\newcommand{\vol}{{\rm {vol}}}
\newcommand{\lct}{{\rm {lct}}}

\newcommand{\Spec}{{\rm {Spec}}}

\newcommand{\dt}{{\rm {dt}}}

\newcommand{\Supp}{{\rm {Supp}}}

\newcommand{\Sym}{{\rm {Sym}}}
\newcommand{\Hom}{{\rm {Hom}}}
\newcommand{\rank}{{\rm {rank}}}
\newcommand{\At}{{\rm {At}}}
\newcommand{\Diff}{{\rm {Diff}}}

\newcommand{\bC}{\mathbb{C}}

\newcommand{\bN}{\mathbb{N}}

\newcommand{\bP}{\mathbb{P}}
\newcommand{\bQ}{\mathbb{Q}}
\newcommand{\bR}{\mathbb{R}}

\newcommand{\mA}{\mathcal{A}}

\newcommand{\mD}{\mathcal{D}}
\newcommand{\mE}{\mathcal{E}}
\newcommand{\mF}{\mathcal{F}}

\newcommand{\mH}{\mathcal{H}}

\newcommand{\mJ}{\mathcal{J}}
\newcommand{\mK}{\mathcal{K}}

\newcommand{\mO}{\mathcal{O}}

\newcommand{\mX}{\mathcal{X}}

\newcommand{\mZ}{\mathcal{Z}}

\newcommand{\tV}{\tilde{V}}

\newcommand{\tg}{\tilde{g}}

\newcommand{\ka}{\mathfrak{a}}

\newcommand{\km}{\mathfrak{m}}
\newcommand{\kn}{\mathfrak{n}}

\begin{document}

\begin{abstract}
In this note, we confirm the first wall $c_1(d)$ of K-moduli of hypersurfaces of degree $d$ in $\mathbb P^n$ for $n\geq 3$ and $d\geq n+1$, where $c_1(d)=\frac{n+1}{nd}$ if $d$ is even; and $c_1(d)=\frac{n+1}{nd-(n+1)}$ if $d$ is odd. Before this wall, the K-moduli is identical to the GIT-moduli for degree $d$ hypersurfaces, so we also call this number the GIT-threshold.  
\end{abstract}

\maketitle
\setcounter{tocdepth}{1}

\tableofcontents

\setcounter{tocdepth}{1}

We work over the field of complex numbers $\bC$.

\section{Introduction}

For hypersurfaces $S_d$ of degree $d$ in $\bP^n$, one can apply either GIT-stability or K-stability to construct moduli spaces. In applying K-stability, we construct moduli for pairs of the form $(\bP^n, cS_d)$ rather than $S_d$ alone, where $c$ is a non-negative rational number such that $(\bP^n, cS_d)$ is a K-semistable log Fano pair. There is a standard wall crossing theory saying that, when $c$ is sufficiently small, the two viewpoints of construction of moduli coincide, i.e., the K-moduli is identical to the GIT-moduli, e.g. \cite[Theorem 1.4]{ADL19}. In this note, we study how small $c$ can be. In other words, we confirm the first wall $c_1(d)$ in this wall crossing picture in the sense of \cite{ADL19, Zhou23}. We will write $c_1$ (resp. $S$) for $c_1(d)$ (resp. $S_d$) since $d$ is clear throughout.

\begin{theorem}\label{thm: main}
Fix $n\geq 3$ and $d\geq n+1$. Define $c_1:=\frac{n+1}{nd}$ if d is even, and $c_1=\frac{n+1}{nd-(n+1)}$ if $d$ is odd. Then a degree d hypersurface $S$ in $\bP^n$ is GIT-(semi/poly)stable if and only if $(\bP^n, cS)$ is K-(semi/poly)stable for any $c\in (0, c_1)$.
\end{theorem}

It is worth mentioning that the case $(n, d)=(2,  d\geq 3)$ is addressed by \cite[Theorem 1.4]{ADL19}, and the case $(n,d)=(3,4)$ is addressed by \cite[Theorem 1.1]{ADL23}. 

We briefly explain the idea of the proof. The central problem is to show that the K-semistable degenerations of pairs $(\bP^n, cS)$ still preserve the ambient space $\bP^n$ for $c\in (0, c_1)$. Let $(\mX, c\mD)\to (0\in C)$ be a $\bQ$-Gorenstein family of K-semistable log Fano pairs over a smooth pointed curve $0\in C$ such that $(\mX_t, \mD_t)$ is of the form $(\bP^n, S)$ for $t\ne 0$.   Put $(X, cD):=(\mX_0, c\mD_0)$. We aim to show that $X\cong \bP^n$. Let $H$ be the hyperplane class on $\bP^n$ and $L$ its limiting class on $X$. We study the rational map $X\dashrightarrow Y\subset \bP^n$ induced by the linear system $|L|$. If $\dim Y=n$, then the rational map induced by $|L|$ is generically finite and the conclusion follows from \cite[Corollary 5.2]{HP24} (see Section \ref{sec: full dim}). If $\dim Y= n-1$ and $d$ is odd, by the analysis on the Cartier index of $L$ (see Section \ref{sec: odd bound}) and the singularities of the base locus $\Sigma$ of $|L|$ (see Section \ref{sec: singularity}), we derive that $L$ is Cartier, which forces $X$ to be isomorphic to $\bP^n$. For the remaining cases, i.e., (1) $\dim Y\leq n-2$ and (2) $\dim Y=n-1$ with $d$ even,  we make use of the slope semistability of the canonical extension of orbifold (co)tangent sheaf (e.g. \cite{XZG26}) to show that they cannot happen (see Section \ref{sec: even case}). 

The notation $X, D, L, Y, \Sigma, n, d, c_1, c$ should be very clear in the above description, and we will employ them throughout the note. 

We acknowledge the use of ChatGPT which found the key idea of lifting the foliation to apply the slope semistability of the canonical extension of the orbifold tangent sheaf of a K-semistable log Fano pair, contained in Section \ref{sec: even case}. The idea in Section \ref{sec: odd case} by applying the singularity property of the linear system $|L|$ (see Section \ref{sec: singularity}) is inspired by some recent AI-assisted works, e.g. \cite{PW26, KSZZ26}. All the contents are written by the authors, who take the full responsibility.

\noindent
\subsection*{Acknowledgement}
CZ is supported by an NSFC grant (No. 12501058) and a grant from Xiamen University (No. X2450214).

\section{Preliminaries}

In this section, we collect some basic knowledge applied in this note. We say $(V, \Delta)$ is a pair if $V$ is a quasi-projective normal variety and $\Delta$ is an effective $\bR$-divisor such that $K_V+\Delta$ is $\bQ$-Cartier. We say a projective pair $(V, \Delta)$ is log Fano if $-K_V-\Delta$ is ample and the pair has klt singularities. 

\subsection{K-stability}

\begin{definition}{\rm (\cite{Fuj19}, \cite{Li17})}
We say a log Fano pair $(V, \Delta)$ is K-semistable if $A_{V, \Delta}(F)-S_{V, \Delta}(F)\geq 0$ for any prime divisor $F$ over $V$. Here,
$$A_{V, \Delta}(F):=\ord_F(K_{\tV}-f^*(K_V+\Delta))+1,$$ 
$$S_{V, \Delta}(F):=\frac{1}{\vol(-K_V-\Delta)}\int_0^\infty\vol(f^*(-K_V-\Delta)-tF){\rm d}t ,$$
where $f: \tV\to V$ is a log resolution such that $F$ is a prime divisor on $\tV$.
\end{definition}

For the K-semistable degeneration $(X, cD)$ as in the introduction, let $F$ be a prime divisor over $X$, then it is not hard to check the following
$$S_{X, cD}(F)=(n+1-cd)\cdot S_L(F). $$

\subsection{Normalized volume}

Given a klt germ $x\in (V, \Delta)$ of dimension $n$, i.e., $(V, \Delta)$ is an $n$-dimensional pair with klt singularities and $x\in V$ is a closed point. The following invariant is introduced by Chi Li (e.g. \cite{Li17}).

\begin{definition}
Let $x\in (V, \Delta)$ be a klt germ as above. Define
\[
\widehat{\vol}(x, V, \Delta)=
\begin{cases}
\text{$\inf_{v} A_{V, \Delta}(v)^n\cdot \vol(v)$} &\quad \text{if $A_{V, \Delta}(v)<\infty$}\\
\text{$\infty$} &\quad \text{if $A_{V, \Delta}=\infty$}\\
\end{cases}
\]
\end{definition}

Readers are referred to \cite{Li17} and \cite{Xu25} for more details. We need several properties of this local invariant.

\begin{theorem}{\rm (\cite{Liu18})}\label{thm: local-global}
Let $(V, \Delta)$ be a K-semistable log Fano pair of dimension $n$ and $x\in V$ is a closed point. Then we have the following local-global volume comparison:
$$\widehat{\vol}(x, V, \Delta)\geq \frac{n^n}{(n+1)^n}\cdot \vol(-K_V-\Delta). $$
\end{theorem}

\begin{theorem}{\rm (\cite{XZ21})}\label{thm: finite degree}
Let $\tilde{x}\in (\tV, \tilde{\Delta})\to x\in (V, \Delta)$ a quasi-\'etale cover of klt germs. Then we have the following finite degree formula:
$$\widehat{\vol}(\tilde{x}, \tV, \tilde{\Delta})=r\cdot \widehat{\vol}(x, V, \Delta), $$
where $r$ is the degree of the quasi-\'etale cover.
\end{theorem}

\begin{theorem}\label{thm: volume}
Let $z\in (Z, cB)$ be an $n$-dimensional klt germ such that $B$ is an effective Cartier divisor and $c$ is a non-negative number. Then we have the following.
\begin{enumerate}
\item $\widehat{\vol}(z, Z, cB)\leq n^n$; if $z\in \Supp\ B$, then $\widehat{\vol}(z, Z, cB)\leq (1-c)n^n$;
\item if $Z$ is locally complete intersection (lci) and singular at $z$, then $\widehat{\vol}(z, Z, cB)\leq 2(n-1)^n$; if moreover $z\in \Supp\ B$, then $\widehat{\vol}(z, Z, cB)\leq 2(n-1-c)^n$.
\end{enumerate}
\end{theorem}

\begin{proof}
These estimates should be well known to experts. We provide a proof for the convenience of the readers.

We first prove (1).  By lower semi-continuity of local volumes (\cite{BL21}), we may assume $z\in Z$ is a smooth point with local coordinates $(x_1,...,x_n)$. Suppose $B$ has a prime component $B_0=(x_1=0)$ with multiplicity $k$ at $z$. Let $v$ be a quasi-monomial valuation centered at $z$ with weight $(\frac{1}{1-ck},1,...,1)$ along the local coordinates $(x_1,...,x_n)$, then 
$$A_{Z, cB}(v)\leq A_{X}(v)-c\cdot v(B_0)=\frac{1}{1-ck}+n-1-\frac{ck}{1-ck}=n,$$
and $\vol(v)=1-ck$. Thus we have
$$\widehat{\vol}(z, Z, cB)\leq A_{Z, cB}(v)^n\cdot \vol(v) \leq n^n\cdot (1-ck)\leq  n^n(1-c).$$

We turn to (2). This is deduced by the same approach as \cite[Theorem 2.12]{Liu22} and we only treat the case $z\in \Supp\ B$ here. Note that $(Z, z)$ is lci and $Z$ is singular at $z$. Let $(Z, z)\to (W,w)$ be the embedding into a smooth $n+r$-fold and denote $\km_w$ the maximal ideal at $w\in W$. Then there exist functions $(f_1,...,f_r)\in \km_w$ locally determining $(Z, z)$.  Let $g$ be the local function of $B$ at $z$ and $\tg$ its lifting at $w$. Choose a general quadric form $Q$ and general independent linear forms $l_2,...,l_r$ and a further general linear form $l_0$ in $\km_w$. We consider the following deformation of $(Z,B, z)$ determined by following equations:
$$f_1+tQ, \ \ \ f_i+tl_i \ (2\leq i\leq r), \ \ \ B_t=(\tg+tl_0=0). $$
For $t\ne 0$, the fiber $(Z_t, B_t, z_t)$ is a nondegenerate quadratic tangent cone, whose normalized volume $\widehat{\vol}(z_t, Z_t, cB_t)$ is no more than $2(n-1-c)^n$. We are done by lower semi-continuity by \cite{BL21}.
\end{proof}

\subsection{Geometry on $X$}

Recall that $(X, cD)$ is a K-semistable degeneration via a $\bQ$-Gorenstein family $(\mX, c\mD)\to C$ of log Fano pairs over a pointed smooth curve $0\in C$, where $(\mX_t, \mD_t)$ is of the form $(\bP^n, S)$ for $t\ne 0$ and $L$ its limiting hyperplane class on $X$. We have the following properties on $X$.

\begin{proposition}\label{prop: geometry}
Notation as above, we have $-K_X\sim (n+1)L$, $D\sim dL$, and $h^0(X, L)=n+1$. In particular, let $X\dashrightarrow Y$ be the rational map induced by $|L|$, then $Y$ is non-degenerate in $\bP^n$. Moreover, if $L$ is Cartier, then $X$ is isomorphic to $\bP^n$.
\end{proposition}

\begin{proof}
Let $\mH^*$ be the hyperplane class on $\mX^*:=\mX\times_{C}(C\setminus 0)$ and $\mH$ the Weil closure of $\mH^*$. So we have $[\mH_0]=[L]$.  Write $\mD^*:=\mD|_{\mX^*}$ and $C^*:=C\setminus 0$. Then we see
$$-K_{\mX^*/C^*}\sim_{C^*} (n+1) \mH^* , \ \ \ \mD^*\sim_{C^*}d\mH^*. $$
Since $X:=\mX_0$ is integral, we easily derive that $-K_X\sim (n+1)L$ and $D\sim dL$. Note that $\mH$ is flat over $C$. By Kawamata-Viehweg vanishing, we have $H^i(X, L)=0$ for $1\leq i\leq n$. Therefore
$$h^0(X, L)=\chi(\mH_0)=\chi(\mH_t)=h^0(\bP^n, \mO_{\bP^n}(1))=n+1. $$
If $L$ is Cartier, then $X$ is isomorphic to $\bP^n$ by \cite{FM21}. The proof is complete.
\end{proof}

\begin{proposition}\label{prop: singular cover}
Suppose $L$ is not Cartier at $x\in X$ and $r$ is the Cartier index of $L$ at $x$. Let $z\in (Z, cB)\to x\in (X, cD)$ be the index one cover associated to $rL$. Then $Z$ is not smooth at $z$.
\end{proposition}

\begin{proof}
We employ the notation in the proof of Proposition \ref{prop: geometry}. Note that $\mX$ has klt singularities, by \cite[Corollary 5.25]{KM98}, $m\mH$ is Cohen-Macaulay for any $m\in \bN$. Thus $\mO_\mX(m\mH)\otimes \mO_X=\mO_X(mL)$ for any $m\in \bN$. In particular, $r\mathcal H$ is also Cartier at $x\in \mathcal \mX$.
Let $x\in U\subset \mX$ be an open neighborhood such that $\mO_U(r\mH)\cong \mO_U$. Consider the following cover
$$\mZ:=\Spec_U\oplus_{i=0}^{r-1}\mO_U(-i\mH)\to U/C. $$
Clearly $\mZ\to U/C$ admits an equivariant $\mu_r$-action over $C$. By the Cohen-Macaulay property mentioned above, $\mZ_0\to U_0$ is exactly identical to the index one cover $z\in Z\to x\in X$, and $\mZ_t\to U_t$ is \'etale $\mu_r$-cover for $t\ne 0$ since $\mH_t$ is already Cartier. Suppose $Z$ is smooth at $z$, then $\mZ\to C$ is smooth near $z\in \mZ$, and $\mZ^{\mu_r}\to C$ is also smooth near $z$, where $\mZ^{\mu_r}$ is the fixed locus under the $\mu_r$-action. Note that $z\in \mZ^{\mu_r}$ since $r$ is the Cartier index of $L$ (resp. $\mH$) at $x$. 
This is a contradiction since $\mZ_t^{\mu_r}$ is empty for $t\ne 0$.
So $Z$ is not smooth at $z$. The proof is complete.
\end{proof}

\subsection{Atiyah sheaf $\At(L)$ on $X$}\label{sec: atiyah lifting}

Readers are referred to \cite[Chapter II, Section 2.2]{Ginzburg98} for more details about this section. Recall that $X$ is a Fano variety and $L$ is a $\bQ$-Cartier Weil divisor on $X$. Replacing $X$ with its regular open subset, we introduce the construction of the Atiyah sheaf $\At(L)$ on $X$. Let $\Diff_{\bC}^1(L, L)$ be the set of $\bC$-linear differential operators on $L$ with order at most one. More precisely, if $s$ is a section of $L$ locally expressed as $s=ge$ under the local frame $e$, then an operator $D\in \Diff_{\bC}^1(L, L)$ acts on $s$ as follows:
$$D(s)=\sigma(D)(g)e+gD(e),$$ 
where $\sigma(D)\in T_X$. 
Write $\At(L):=\Diff_{\bC}^1(L, L)$ and there exists the natural projection:
$$\sigma: \At(L)\to T_X, \ \ \ D\mapsto \sigma(D).$$
Suppose $\sigma(D)=0$, then $D$ is $\mO_X$-linear and thus contained in $\mH om_{\mO_X}(L, L)\cong \mO_X$. So we have the following short exact sequence:
$$0\to \mO_X\to \At(L)\to T_X\to 0. $$
Let $U_i\subset X$ be open subsets such that $L|_{U_i}\cong \mO_{U_i}$, then locally on $U_i$ we have
$$\At(L)|_{U_i}\cong \mO_{U_i}\oplus T_{U_i}. $$
Suppose $e_i$ is the frame on $U_i$ and $s|_{U_i}=f_ie_i$, where $e_j=g_{ij}e_i$ and $\{g_{ij}\}$ are the transition functions associated to $L$. Then $D|_{U_i}$ can be expressed in the form $(a_i, v_i)$, where $a_i\in \mO_{U_i}$, $v_i\in T_{U_i}$, and
$$(a_i, v_i)(f_ie_i)=\bigg(a_if_i+v_i(f_i)\bigg)e_i. $$
On $U_i\cap U_j$, we also have
\begin{align*}
(a_j, v_j)(f_je_j)&=\bigg(a_jf_j+v_j(f_j)\bigg)e_j \\
&=\bigg(a_jf_i+v_j(\frac{f_i}{g_{ij}})\cdot g_{ij}\bigg)e_i\\
&=\bigg(a_jf_i+v_j(f_i)-\frac{v_j(g_{ij})}{g_{ij}}\cdot f_i\bigg)e_i.
\end{align*}
This means on $U_i\cap U_j$ we have
$$(a_i, v_i)=(a_j-\frac{v_j(g_{ij})}{g_{ij}}, v_j) .$$
Therefore, $\At(L)$ can also be viewed as the gluing of $\{\mO_{U_i}\oplus T_{U_i}\}$ with respect to the transition functions $\{\frac{{\rm d} g_{ij}}{g_{ij}}\}$. So the extension $\At(L)$ is obtained under the extension class in $\bQ\cdot L$.

\section{Index bound for odd degrees}\label{sec: odd bound}

Let $r$ be the Cartier index of $L$ at $x\in (X, cD)$. Let $(Z, cB;z)\to (X, cD; x)$ be the index one cover associated to $rL$. In this section we assume $n\geq 3$ and $d\geq n+1$, but we only focus on odd $d$. Also note that $c\in (0, \frac{n+1}{nd-(n+1)})$ for odd $d$.

\begin{proposition}\label{prop: odd bound}
Suppose $n\geq 3$ and $d\geq n+1$. Assume $d$ is odd and $(n,d)\ne (4,5)$. Then $r\leq 3$.
\end{proposition}

\begin{proof}
By Theorem \ref{thm: local-global} we have
$$\widehat{\vol}(x, X, cD)\geq (\frac{n}{n+1})^n\cdot (n+1-cd)^n. $$
If $x\in \Supp\ D$, then applying Theorem \ref{thm: volume},
$$\widehat{\vol}(z, Z, cB)\leq (1-c)n^n. $$
By finite degree formula, i.e., Theorem \ref{thm: finite degree}, we have the following estimate for $d\geq 5$:
$$r\leq \frac{(1-c)n^n}{(\frac{n}{n+1})^n\cdot (n+1-cd)^n}<4,$$
implying $r\leq 3$ for this case.

If $x\notin\Supp\ D$, then $r\mid d$ since $D\sim dL$. We first assume $d\geq 7$ and now we have
$$\widehat{\vol}(z, Z, cB)\leq n^n ,$$
and thus we see the following estimate for $d\geq 7$:
$$r\leq \frac{n^n}{(\frac{n}{n+1})^n\cdot (n+1-cd)^n}<5. $$
Since $d$ is odd and $r\mid d$ in this case, we see $r\leq 3$. Next we assume $d=5$, for which we only need to consider $(n,d)=(4,5)$ and $(n, d)=(3,5)$. For $(n,d)=(4,5)$, we have the following estimate
$$r\leq \frac{n^n}{(\frac{n}{n+1})^n\cdot (n+1-cd)^n}<6. $$
So it may happen $r=5$ for $(n, d)=(4,5)$. 
However, since $(n,d)\ne (4,5)$ by assumption, the only remaining case is $(n,d)=(3,5)$. In dimension three, we have better estimate by \cite{LX19}:
$$\widehat{\vol}(z, Z, cB)\leq 2(n-1)^n,$$
which implies the following for $(n,d)=(3,5)$:
$$r\leq \frac{2(n-1)^n}{(\frac{n}{n+1})^n\cdot (n+1-cd)^n}<4. $$
Above all, we have $r\leq 3$ and the proof is complete.
\end{proof}

\begin{proposition}\label{prop: lci odd bound}
Suppose $n\geq 3$ and $d\geq n+1$. Assume $d$ is odd, and $Z$ is singular and lci at $z$. Then $r= 1$.
\end{proposition}

\begin{proof}
The proof is similar to that of Proposition \ref{prop: odd bound}. First we have
$$\widehat{\vol}(x, X, cD)\geq (\frac{n}{n+1})^n\cdot (n+1-cd)^n. $$
If $x\in \Supp\ D$, then 
$$\widehat{\vol}(z, Z, cB)\leq 2(n-1-c)^n. $$
By finite degree formula, we have the following estimate for $d\geq 5$:
$$r\leq \frac{2(n-1-c)^n}{(\frac{n}{n+1})^n\cdot (n+1-cd)^n}<2,$$
implying $r=1$ in this case.

If $x\notin\Supp\ D$ and $r\ne 1$, then $r\mid d$ since $D\sim dL$. Now we have
$$\widehat{\vol}(z, Z, cB)\leq 2(n-1)^n, $$
and thus we see the following estimate for $d\geq 5$:
$$r\leq \frac{2(n-1)^n}{(\frac{n}{n+1})^n\cdot (n+1-cd)^n}<4. $$
This forces $r=2$ or $r=3$. But $r\mid d$ and $d$ is odd, we see $r\ne 2$; if $r=3$, then $d\geq 9$ and  
$$r\leq \frac{2(n-1)^n}{(\frac{n}{n+1})^n\cdot (n+1-cd)^n}<3,$$
a contradiction. Thus $r=1$. 
\end{proof}

\section{Singularities of the linear system $|L|$}\label{sec: singularity}

Consider the following diagram:

\begin{center}
	\begin{tikzcd}[column sep = 2em, row sep = 2em]
	&W\arrow[ld,"p",swap]\arrow[rd,"q"]&\\
	 X \arrow[rr,"|L|",dashed]&& Y\subset \bP^n,
        \end{tikzcd}
\end{center}
where $p: W\to X$ is a resolution such that $P:=p^*L=M+E$ is the sum of free part and fixed part, and $q: W\to Y$ is induced by $M$. Put $a_i:=P^{n-i}M^i$ for $0\leq i\leq n$. We easily see 
$$1=P^n=a_0\geq a_1\geq\cdots\geq a_{n}=M^n\geq 0. $$
We have the following characterization of singularities of base locus of $|L|$. 

\begin{proposition}\label{prop: singularities}
Let $F$ be a prime divisor over $X$ centered in $p(E)$, then 
$$A_X(F)\geq \kappa_c\cdot \ord_F(E), $$ 
where $\kappa_c:=(n+1-cd)\cdot \frac{\sum_i a_i}{n+1}$.
\end{proposition}

\begin{proof}
We may assume $F$ is a prime divisor on $W$ up to resolving $W$. Denote $\beta:=\ord_F(E)>0$ (note that $F$ is centered in $p(E)$). 

We have the following computation:
\begin{align*}
S_L(F)&=\frac{1}{L^n}\int_0^\infty \vol(P-tF)\dt\\
&=\int_0^\infty \vol((1-\frac{t}{\beta})P+\frac{t}{\beta}M+\frac{t}{\beta}E-tF)\dt\\
&\geq \beta\int_0^1 \vol((1-u)P+uM){\rm d}u\\
&=\beta\int_0^1 ((1-u)P+uM)^n {\rm d}u\\
&=\frac{\sum_{i=0}^n a_i}{n+1}\cdot \beta.
\end{align*}
By the K-semistability of $(X, cD)$, we have
$$A_{X, cD}(F)=A_X(F)-c\cdot \ord_F(D)\geq S_{X, cD}(F)=(n+1-cd)\cdot S_L(F). $$
Thus we see
$$A_X(F)\geq (n+1-cd)\cdot\frac{\sum_i a_i}{n+1}\cdot \ord_F(E). $$
Putting $\kappa_c:=(n+1-cd)\cdot \frac{\sum_i a_i}{n+1}$, we have $A_X(F)\geq \kappa_c\cdot \ord_F(E)$. The proof is complete. 
\end{proof}

The formula in Proposition \ref{prop: singularities} means $X$ has controlled singularities along the base locus $p(E)$. In fact, for any $x\in p(E)$, suppose the Cartier index of $L$ is $r$ at $x$. Let $\pi: (Z, cB; z)\to (X, cD; x)$ be the index-one cover associated to $rL$. Denote
$$I:={\rm Im}\  (\Sym^r H^0(X, L)\otimes \mO_X(-rL)\to \mO_X), $$
then $\ord_F(I)=r\cdot \ord_F(E)$ for any prime divisor $F$ over $X$ centered in $p(E)$. Thus $(X, I^\frac{\kappa_c}{r})$ is lc along $p(E)$. Put  $J^r:=I\mO_Z$. Note that $(Z, cB)\to (X, cD)$ is a quasi-\'etale cover, then we see $(Z, J^{\kappa_c})$ is lc.

\begin{proposition}\label{prop: lci}
Notation as above. We have the following conclusions.
\begin{enumerate}
\item Let $s$ be the codimension of $p(E)$ in $X$, then $s\geq \kappa_c$;
\item If $\kappa_c>n-2$, then $Z$ is locally complete intersection at $z$; 
\end{enumerate}
\end{proposition}

\begin{proof}
Let $\eta$ be the generic point of a component of $p(E)$ of codimension $s$, and let $\xi\in Z$ be the generic point over $\eta$. Locally we write $(R, \km):=(Z, \xi)$. Since $Z$ has canonical Gorenstein singularities, $(R, \km)$ is a canonical Gorenstein local ring of dimension $s$.  First observe that $J\subset \km$. Thus we see $\kappa_c\leq \lct(J)\leq \lct(\km)$. By the lower semi-continuity of log canonical thresholds we know $\lct(\km)\leq s$, hence $\kappa_c\leq s$. This is (1). 

Now we turn to (2) and locally write $(R', \km'):=(Z, z)$.  Since $J\subset \km'$, we still have $\kappa_c\leq \lct(J)\leq \lct(\km')$. Assume $\kappa_c>n-2$, then $\lct(\km')>n-2$. Let $\ka=(f_1,...,f_n)$ be a reduction of $\km'$, and $\epsilon$ is a small positive number such that $\lct(\km')>n-2+\epsilon$. Then we first see
$$\km'^2=\km'^2\mJ(\km'^{n-2+\epsilon})\subset \mJ(\km'^{n+\epsilon}).$$
On the other hand, since $\ka$ is a reduction of $\km'$, we have the following by \cite[Corollary 1.4]{LLS08}
$$\mJ(\km'^{n+\epsilon})=\ka\mJ(\km'^{n-1+\epsilon})\subset \ka. $$
Above all we have $\km'^2\subset \ka$. Let $\kn:=\km'/\ka$ be the maximal ideal of $\bar{R'}:=R'/\ka$. If $\kn\neq 0$, then $\bar{R'}$ is a Gorenstein Artin local ring and $\kn^2=0$. Thus we see 
$$\Hom_{\bar{R'}}(\bar{R'}/\kn, \bar{R'})=(0:\kn)=\kn.$$
By the characterization of a Gorenstein Artin local ring (e.g. \cite[Theorem 18.1(4)]{Matsumura89}, we know that $\dim_{R'/\km'}(\kn)= 1$, which implies that $\km'$ is generated by at most $n+1$ elements. Thus the embedding dimension of $(R', \km')$ is at most $n+1$. Since $R'$ is Cohen-Macaulay, we see that $(R', \km')$ is either regular or a hypersurface. In particular, $(R', \km')$ is lci. The proof is complete.
\end{proof}

\section{The first case: $\dim Y= n$}\label{sec: full dim}

\begin{theorem}\label{thm: full dim}
Assume $\dim Y=n$. Then $L$ is Cartier and $X$ is isomorphic to $\bP^n$.
\end{theorem}

\begin{proof}
 Implied by \cite[Corollary 5.2]{HP24}.
\end{proof}

\section{The second case: $\dim Y=n-1$ and $d$ is odd}\label{sec: odd case}

In this section, we assume $n\geq 3$ and $d\geq n+1$, and we also assume $\dim Y=n-1$ and $d$ is odd. By the analysis in Section \ref{sec: odd bound}, the Cartier index of $L$ is at most $3$ except the case $(n,d)=(4,5)$ (see Proposition \ref{prop: odd bound}). We first have the following easy property for the intersection numbers $a_{n-1}:=PM^{n-1}$ (see Section \ref{sec: singularity}).

\begin{lemma}\label{lem: a}
Assume $n\geq 3$, $d\geq n+1$, $\dim Y=n-1$, $d$ is odd, and $(n,d)\ne (4,5)$. Then $a_{n-1}\geq \frac{2}{3}$.
\end{lemma}

\begin{proof}
Recall the diagram in Section \ref{sec: singularity}. By Proposition \ref{prop: odd bound} we know that the Cartier index of $L$ is at most $3$. Thus the coefficients of $P=M+E$ (see the notation in Section \ref{sec: singularity}) are at least $\frac{1}{3}$. Since $L$ is ample, some component of $E$ must dominate $Y$, which implies $q_* 3P\geq [Y]$. So we have the following by projection formula
$$3PM^{n-1}=3P(q^*(\mO_Y(1)))^{n-1}=(q_*3P)(\mO_Y(1))^{n-1}\geq \deg Y\geq 2. $$
The last inequality holds since $Y$ is non-degenerate in $\bP^n$. Hence we have 
$$a_{n-1}\geq  \frac{2}{3}.$$
The proof is complete.
\end{proof}

The above lower bound for $a_{n-1}$ will lead to a lower bound for $\kappa_c:=\frac{(n+1-cd)\sum_{i}a_i}{n+1}$, since $1=a_0\geq a_1\geq\cdots\geq a_{n-1}$. Below we give a finer bound. 

\begin{lemma}\label{lem: lci lemma}
Assume $n\geq 3$, $d\geq n+1$, $\dim Y=n-1$, and $d$ is odd. Then $\kappa_c>n-2$.
\end{lemma}

\begin{proof}
Let $s$ be the codimension of $p(E)$ in $X$ (see the notation in Section \ref{sec: singularity}). If $s>k$, then $P^{n-l}E^l=0$ for $1\leq l\leq k$ by projection formula. Thus we easily derive
$$P^n=a_0=a_1=\cdots =a_k=1. $$
Assume $(n,d)\neq (4,5)$. By Lemma \ref{lem: a}, we already have
$$a_{k+1}\geq\cdots a_{n-1}\geq \frac{2}{3}. $$
For any $c\in (0, c_1)$, we have the following estimate for $\kappa_c$:
$$\kappa_c\geq (1-\frac{cd}{n+1})\cdot \bigg((k+1)+\frac{2}{3}(n-(k+1))\bigg). $$
It is not hard to verify $\kappa_c>k+1$ for $0\leq k\leq n-3$. By Proposition \ref{prop: lci}, we have $s\geq \kappa_c>k+1$. Thus we see $\kappa_c>n-2$.

Assume $(n, d)=(4,5)$. By the proof of Proposition \ref{prop: odd bound}, $r\leq 5$. By the same proof of Lemma \ref{lem: a}, we know 
$$a_{n-1}=a_3\geq \frac{2}{5}.$$
We hence have the following:
$$a_0=1,\ \ \ a_1\geq a_2\geq a_3\geq \frac{2}{5},$$
which implies the following estimate for $(n,d)=(4,5)$:
$$\kappa_c\geq (1-\frac{cd}{n+1})\sum_ia_i\geq (1-c)(1+\frac{2}{5}+\frac{2}{5}+\frac{2}{5})>n-3=1. $$
By Proposition \ref{prop: lci}, the codimension of $p(E)$ is bigger than $1$, and we thus have $a_0=a_1=1$. 
By the convexity of the intersection numbers (e.g. \cite[Example 1.6.4]{Laz04a}), we know
$$a_2^2\geq a_1\cdot a_3\geq \frac{2}{5}, $$
which implies $a_2>\frac{3}{5}$ and the following finer estimate for $(n,d)=(4,5)$:
$$\kappa_c\geq (1-\frac{cd}{n+1})\sum_ia_i> (1-c)(1+1+\frac{3}{5}+\frac{2}{5})>n-2=2. $$
The proof is complete.
\end{proof}

\begin{theorem}\label{thm: odd thm}
Assume $n\geq 3$, $d\geq n+1$, and $d$ is odd. Then it cannot happen that $\dim Y=n-1$. 
\end{theorem}

\begin{proof}
Assume $\dim Y=n-1$ and $d$ is odd. Suppose $L$ is not Cartier at $x\in X$ and let $r$ be the Cartier index of $L$ at $x$. Let $z\in (Z, cB)\to x\in (X, cD)$ be the index one cover associated to $rL$. By Lemma \ref{lem: lci lemma} and Proposition \ref{prop: lci}, $Z$ is lci at $z$. By Proposition \ref{prop: singular cover}, $Z$ is not smooth at $z$. Applying Proposition \ref{prop: lci odd bound}, we know that $L$ is Cartier at $x$, leading to a contradiction by Proposition \ref{prop: geometry}.
\end{proof}

\section{The remaining cases}\label{sec: even case}

Since we focus on $n\geq 3$ and $d\geq n+1$ in this note, the remaining cases are (1) $\dim Y\leq n-2$; and (2) $\dim Y=n-1$ with even $d$. We show in this section that these cases cannot happen. Write $m:=\dim Y$ and assume $1\leq m\leq n-1$ in this section. Consider the following diagram:
\begin{center}
	\begin{tikzcd}[column sep = 2em, row sep = 2em]
	 \Gamma \arrow[d,"p_1", swap] \arrow[rr,"q_1"]&& Y'\arrow[d,"\phi"]&\\
	 X\arrow[rr,"|L|", dashed]&& Y&\subset \bP^n,
	\end{tikzcd}
\end{center}
where $\phi: Y'\to Y$ is a smooth resolution, $\Gamma$ is a sufficiently high smooth birational model of $X$ such that $p_1^*L:=M_1+E_1$ is the sum of free part and fixed part, and $q_1$ is induced by $|M_1|$. Fix the following notation:
$$B_0:=p_1(E_1), \ \ M_0:=p_1(M_1),\ \ A:=\phi^*{\mO_Y(1)}.  $$
Here $B_0$ and $M_0$ are respectively the fixed divisor and the movable divisor of $|L|$.

Let $\mE$ be a torsion free coherent sheaf on $X$, we define the degree and slope of $\mE$ with respect to $L$ as follows:
$$\deg_L(\mE):= c_1(\mE)\cdot L^{n-1},\ \ \ \mu_L(\mE):=\frac{\deg_L(\mE)}{\rank\ \mE}=\frac{c_1(\mE)\cdot L^{n-1}}{\rank\ \mE} .$$
Let $\mF$ be the algebraic foliation induced by the rational map $X\dashrightarrow Y$. We will derive two incompatible estimates of the degree of $\mF$ to show that the two remaining cases cannot appear, as in the following theorem.

\begin{theorem}\label{thm: even thm}
Assume $n\geq 3$, $d\geq n+1$. Then the following two cases cannot happen: (1) $\dim Y=n-1$ and $d$ is even; (2) $\dim Y\leq n-2$.
\end{theorem}

\subsection{Lower bound of $\deg_L(\mF)$}

 We first prove the following independent lemma.

\begin{lemma}\label{lem: global generation}
Let $T$ be a smooth projective variety of dimension $m$ and $G$ a globally generated big line bundle on $T$ such that $h^0(T, G)\geq m+2$. Then $K_T+mG$ is effective. 
\end{lemma}

\begin{proof}
The lemma is clear when $m=1$, so we assume $m\geq 2$. Since $|G|$ is globally generated, we successively use general members of $|G|$ to cut $T$ for $m$ times and obtain the following flag:
$$T_m\subset T_{m-1}\subset\cdots\subset T_1\subset T_0=T,$$
where $T_{m-1}$ is a smooth curve and $T_m$ is a dimension zero reduced subset. Write $C:=T_{m-1}$. By successively applying Kawamata-Viehweg vanishing and adjunction, we have the following surjective map:
$$ H^0(T, K_T+mG)\to H^0(C, K_C+G|_C). $$
On the other hand, we have the following 
$$m+2\leq h^0(T, G)\leq h^0(C, G|_C)+m-1, $$
which implies $h^0(C, G|_C)\geq 3$, and thus $h^0(C, K_C+G|_C)>0$.
This implies that $K_T+mG$ is effective by the above surjective map. The proof is complete.
\end{proof}

\begin{lemma}\label{lem: lower bound}
Assume $m:=\dim Y \leq n-1$. Then we have 
$$\deg_L(\mF)\geq n+1-m+m\deg_L(B_0).$$
\end{lemma}

\begin{proof}
Denote $\mF'$ to be the foliation associated to $q_1: \Gamma\to Y'$. We clearly have
$$c_1(\mF')=-K_\Gamma+q_1^*K_{Y'}+R $$
for some effective divisor $R$. Thus 
$$c_1(\mF)=p_{1*}(c_1(\mF'))\geq -K_X+p_{1*}q_1^*K_{Y'}. $$
Notice that $Y'$ is smooth of dimension $m$ and $A$ is big and globally generated on $Y'$. By Lemma \ref{lem: global generation}, we see 
$$K_{Y'}\geq -mA. $$
We thus have the following
\begin{align*}
c_1(\mF)&\geq -K_X-m\cdot p_{1*}q_1^*A\\
&=-K_X-m\cdot M_0\\
&=-K_X-m(L-B_0) \\
&=(n+1-m)L+mB_0.
\end{align*}
This directly implies $\deg_L(\mF)\geq n+1-m+m\deg_L(B_0)$ as desired.
\end{proof}

\subsection{Upper bound of $\deg_L(\mF)$}

Remember that $(X, cD)$ is a K-semistable log Fano pair. By \cite[Section 5]{XZG26}, there is a canonical extension $\mE_{X, cD}$ as follows
$$0\to \Omega_X(\log (cD))\to \mE_{X, cD}\to \mO_X\to 0 $$
associated to the extension class in $\bQ\cdot L$, where $\Omega_X(\log (cD))$ is the orbifold cotangent sheaf with respect to $(X, cD)$. Dualizing the above sequence, we write $\mA_{X, cD}:=\mE_{X, cD}^{\lor}$. We then have the following
$$0\to \mO_X\to \mA_{X, cD}\to T_{X, cD}\to 0, $$
where $T_{X, cD}$ is the orbifold tangent sheaf. The above sheaves and exact sequences should be understood on the regular locus of $X$, forming a big open subset. By \cite[Theorem 5.8 and Appendix B]{XZG26}, we have the following result:

\begin{lemma}\label{lem: slope ss}
The canonical extension $\mA_{X, cD}$ of the orbifold tangent sheaf $T_{X, cD}$ is slope semistable with respect to $L$.
\end{lemma}

Notice that $c_1(\mA_{X, cD})=-K_X-cD$ and $\rank (\mA_{X, cD})=n+1$. Thus we have 
$$\deg_L(\mA_{X, cD})=n+1-cd, \ \ \ \mu_L(\mA_{X, cD})=\frac{n+1-cd}{n+1}.$$
To make use of the slope semistability of $\mA_{X, cD}$ to control $\deg_L(\mF)$, we first need to lift $\mF$ to a subsheaf $\tilde{\mF}\subset \mA_{X, cD}$, and then compare $\deg_L(\mF)$ and $\deg_L(\tilde{\mF})$. We first construct $\tilde{\mF}\subset \mA_{X, cD}$. Write
$$cD=\sum_P \delta_P P, \ \ \ \ \  B_0:=\sum_P b_P P, $$
where $P$ runs over prime divisors on $X$. It is clear that $0\leq \delta_P<1$ and $b_P\in \bN$.
By the construction in Section \ref{sec: atiyah lifting}, there is an extension of $T_X$ with extension class in $\bQ\cdot L$ as follows:
$$0\to \mO_X\to \At(L)\to T_X\to 0. $$
We first construct the lifting of $\mF\subset T_X$ to $\At(L)\otimes_{\mO_X} \mK_X$, where $\mK_X$ is the sheaf of rational functions on $X$. Suppose the prime divisor $P$ is generically locally defined by $x=0$, and $s_0$ is a section of $L$ that can be locally written as $s_0=g_0e_0$ under the local frame $e_0$. If $v\in T_X$, then its lifting $P_v\in \At(L)$ acts on the section $s_0$ locally as
$$P_v(s_0)=v(g_0)e_0+g_0P_v(e_0). $$
From now on we assume $v\in \mF$, then its lifting $P_v$ satisfies $P_v(s_0)=0$. That is to say,
$$v(g_0)e_0+g_0P_v(e_0)=0, $$
which gives
$$P_v(e_0)=-\frac{v(g_0)}{g_0}e_0. $$
For any other section $s$ of $L$, it is direct to check
$$P_v(s)=P_v(\frac{s}{s_0}\cdot s_0)= v(\frac{s}{s_0})s_0+\frac{s}{s_0}\cdot P_v(s_0)=0,$$
since $v(\frac{s}{s_0})=0$.
Therefore, the lifting $P_v$ locally has the form
$$(\clubsuit)\ \ \ \ (-\frac{v(g_0)}{g_0}, v)\in \mK_U\oplus T_U \ \ \ \ \text{(see Section \ref{sec: atiyah lifting})}$$
for some open subset $U$ of the generic point of $P$, where $\mK_U$ is the sheaf of rational functions on $U$. 
So $P_v$ is not necessarily contained in $\mA_{X, cD}$, since the corresponding element in $\mA_{X, cD}$ is locally of the form: 
$$(\clubsuit\clubsuit)\ \ \ \ (\text{a regular function}, x^{\delta_P}\partial_x+\text{other directions}).$$ 
We aim to compare $(\clubsuit)$ and $(\clubsuit\clubsuit)$. If $P$ is contained in $\Supp\ B_0$, i.e., $b_P>0$, then $g_0=x^{b_P}u$ for some unit $u$ at the generic point of $P$, and thus
$$-\frac{v(g_0)}{g_0}= - \bigg(b_P\cdot\frac{v(x)}{x}+\frac{v(u)}{u}\bigg);$$ 
if $P$ is not contained in $\Supp\ B_0$, i.e., $b_P=0$, then we see
$$-\frac{v(g_0)}{g_0}= - \frac{v(u)}{u}. $$
At the generic point of the prime divisor $P$, we may choose a local basis $\{v_i\}_{i=1}^{\rank\ \mF}$ of $\mF$ such that $v_i(x)=0$ for $i\geq 2$. Thus, only the lifting of $v_1$ may require the multiplication by $x$ or $x^{\delta_P}$, and the determinant of $\tilde{\mF}$ loses at most one such factor along $P$.
Inspired by this, we define two auxiliary divisors as follows: 
$$\Theta_1:=\sum_{P, \ \ b_P>0} P, \ \  \text{and}\ \  \Theta_2:=\sum_{P, \ \ b_P=0}\delta_P P.$$
We automatically have $\Theta_1\leq B_0$ and $\Theta_2\leq cD$.

\begin{lemma}\label{lem: upper bound}
Notation as above, we have
$$\deg_L(\tilde{\mF})\geq \deg_L(\mF)-\deg_L(B_0+cD)$$
and
$$\deg_L(\tilde{\mF})\leq \frac{n-m}{n+1}\cdot \deg_L(\mA_{X, cD})= \frac{n-m}{n+1}\cdot (n+1-cd).$$
In particular, 
$$\deg_L(\mF)\leq \frac{n-m}{n+1}\cdot (n+1-cd)+\deg_L(B_0)+cd. $$
\end{lemma}

\begin{proof}
By the construction of $\tilde{\mF}$, we see the following
$$c_1(\tilde{\mF})\geq c_1(\mF)-(\Theta_1+\Theta_2) \geq c_1(\mF)-(B_0+cD),$$
from which we conclude the first inequality. The second inequality follows from the slope semistability of $\mA_{X, cD}$ (e.g. Lemma \ref{lem: slope ss}) and the fact that $\rank\ \mF=\rank\ \tilde{\mF}$:
$$\frac{\deg_L(\tilde{\mF})}{\rank\ \tilde{\mF}}\leq \frac{\deg_L(\mA_{X, cD})}{\rank\ \mA_{X, cD}}. $$
The third inequality follows directly from the combination of the first and the second inequalities.
\end{proof}

\begin{proof}[Proof of Theorem \ref{thm: even thm}]
By the two estimates in Lemma \ref{lem: lower bound} and Lemma \ref{lem: upper bound}, we have the following:
$$n+1-m+m\deg_L(B_0)\leq \frac{n-m}{n+1}\cdot (n+1-cd)+\deg_L(B_0)+cd. $$
Equivalently, we have the following:
$$ (\spadesuit)\ \ \ 1+(m-1)\deg_L(B_0)\leq cd\cdot (1-\frac{n-m}{n+1}).$$
Since $m\leq n-1$, if $d$ is even, then $c\in (0, \frac{n+1}{nd})$ and
$$cd\cdot (1-\frac{n-m}{n+1})< \frac{n+1}{n}\cdot \frac{1+m}{n+1}\leq 1, $$
which contradicts to $(\spadesuit)$. If $d$ is odd and $m\leq n-2$, then $c\in (0, \frac{n+1}{nd-(n+1)})$ and
$$cd\cdot (1-\frac{n-m}{n+1})< \frac{(n+1)d}{nd-(n+1)}\cdot \frac{1+m}{n+1}\leq 1, $$
which also contradicts to $(\spadesuit)$. This says the two cases in Theorem \ref{thm: even thm} cannot appear. The proof is complete.
\end{proof}

\section{Optimal bound}

Combining Theorem \ref{thm: full dim}, Theorem \ref{thm: odd thm}, and Theorem \ref{thm: even thm}, we get the following result at once.

\begin{theorem}\label{thm: rigidity}
Assume $n\geq 3$ and $d\geq n+1$. Define $c_1:=\frac{n+1}{nd}$ if d is even, and $c_1=\frac{n+1}{nd-(n+1)}$ if $d$ is odd. Fix any rational $c\in (0, c_1)$. Let $(\mX, c\mD)\to (0\in C)$ be a $\bQ$-Gorenstein family of K-semistable log Fano pairs over a smooth pointed curve $0\in C$ such that $(\mX_t, \mD_t)$ is of the form $(\bP^n, S_d)$ for $t\ne 0$. Then the special degeneration $\mX_0$ is isomorphic to $\bP^n$ and $\mD_0$ is a degree $d$ hypersurface. 
\end{theorem}

We next explain that $c_1$ in Theorem \ref{thm: rigidity} is indeed optimal. 

\begin{theorem}\label{thm: optimal}
Assume $n\geq 3$ and $d\geq n+1$. Notation as in Theorem \ref{thm: rigidity} and write $(X, D):=(\mX_0, \mD_0)$. If $c=c_1$, then it can happen that the special degeneration $X$ is not isomorphic to $\bP^n$. 
\end{theorem}

\begin{proof}
This should be well-known to experts. It suffices to find a K-semistable degeneration at $c_1$ such that the special degeneration is not $\bP^n$.
Define
  
\begin{equation*}
\Delta_d=
\begin{cases}
\frac{d}{2}Q& \text{$d$ even}\\
\frac{d-1}{2}Q+L& \text{$d$ odd}\\
\end{cases}
\end{equation*}
where $Q\subset \bP^n$ is a smooth quadric hypersurface and $L\subset \bP^n$ is a hyperplane such that $(\bP^n, Q+L)$ is log smooth. By \cite[Theorem 1.4]{LZ23}, one can easily check that $(\bP^n, c_1\Delta_d)$ is K-semistable, and further we have
$$A_{X, c_1\Delta_d}(Q)-S_{X, c_1\Delta_d}(Q) =0.$$
By \cite{LWX21}, $Q$ induces a special test configuration whose central fiber is also K-semistable. However, the central fiber, which is a cone over the quadric hypersurface, is clearly not isomorphic to $\bP^n$.
\end{proof}

\begin{proof}[Proof of Theorem \ref{thm: main}]
After establishing Theorem \ref{thm: rigidity} and Theorem \ref{thm: optimal}, Theorem \ref{thm: main} will be derived in a standard way, e.g. \cite[Proof of Theorem 9.24]{ADL19}.
\end{proof}

\begin{remark}
Since we mainly concern $d\geq n+1$ in this note, one may naturally ask whether this explicit GIT-threshold still applies to lower degree case. This is however not true. For example, taking $(n,d)=(4, 3)$, then $\frac{n+1}{nd-(n+1)}=\frac{5}{7}$.  But the GIT-threshold for cubic threefolds is $\frac{5}{6}$, e.g. \cite[Proof of Proposition 5.2]{SZZ26}. It is an interesting question to further compute the GIT-thresholds for lower degree hypersurfaces. 
\end{remark}

\bibliography{reference.bib}
\end{document}